\documentclass{amsart}

\usepackage{fancyhdr}
\newtheorem{theorem}{Theorem}[section]
\newtheorem{lemma}[theorem]{Lemma}

\theoremstyle{definition}

\theoremstyle{remark}
\newtheorem{remark}[theorem]{Remark}

\numberwithin{equation}{section}

\begin{document}

\title{Complex Monge-Amp\`ere Equations with Oblique Boundary Value}

\author{Ni Xiang }
\address{Faculty of Mathematics $\&$ Computer Science, Hu Bei University,
P.R. China}
\email{nixiang@hubu.edu.cn}

\author{Xiao-Ping Yang}
\address{School of Science, Nanjing University of Sci. $\&$ Tech., Nanjing, P.R.China 210094}
\email{xpyang@mail.njust.edu.au}
\thanks{The first author was
supported by the National Natural Science Foundation of China
(No.11101132) and foundation of Hubei Provincial department of
education(No.Q20120105).}





\begin{abstract}
The existence and regularity of the classical plurisubharmonic
solution for complex Monge-Amp\`ere equations subject to the
semilinear oblique boundary condition which is $C^1$ perturbation of
the Neumann boundary condition, are proved in the certain strictly
pseudoconvex domain in $\bf C^n$.
\end{abstract}

\maketitle
\section{ Introduction and main results}
Let $\Omega$ be a bounded strictly pseudoconvex domain in $\bf C^n$
with $C^\infty$ boundary, $f(z)$ be a positive function on $\Omega$,
and $\varphi(z,u)$ be a function on $\partial \Omega\times \bf R$.
We shall study the existence and regularity of plurisubharmonic
solutions to the complex Monge-Amp\`ere equations:
\begin{equation}\label{cma}
\det\frac{\partial^2u}{\partial z_i\partial\overline{z_j}}=f(z), \ \
\ \ z\in\Omega, \end{equation}
 with the oblique boundary condition:
\begin{equation}\label{ob1}
D_\beta u=\varphi(z,u),\ \ \ \ \  z\in
\partial\Omega . \end{equation}

Assume
 $\beta$ is a strictly oblique unit
 vector field, satisfying  \begin{equation}\label{ob2'}\beta(z)\cdot\nu(z)\geq\beta_0>0,\ \
 \end{equation} where $ \nu(z) $ is the unit outer normal vector
field to $\partial\Omega$. In order to apply the method of
continuity for the existence of solutions, it is necessary to obtain
a priori estimates. We also need to assume some conditions on
$\varphi$ to derive the maximum modulus estimates. We assume the
conditions in \cite{LTU,U2} namely
\begin{equation}\label{ob4}\varphi(z,u)<0\ for\ all \ z\in \partial\Omega\ \ \  and \ \  all \ \ u\geq N_1, \end{equation}
for some constant $N_1$, and
\begin{equation}\label{ob5}\varphi(z,u)\rightarrow \infty \ \ as \ \
u\rightarrow-\infty \ \ uniformly \ \ for \ \ z\in
\partial\Omega.
\end{equation}Conditions (\ref{ob4}) and
(\ref{ob5}) are used for the upper bound and the lower bound of the
solution, respectively. We also assume that
\begin{equation}\label{ob6} \inf
f(\overline{\Omega})>0.\end{equation}

The Dirichlet problem for complex Monge-Amp\`ere equations
(\ref{cma}) has been an object of intensive research. In 1976,
Bedford and Taylor in \cite{BT} considered weak plurisubharmonic
solutions for complex Monge-Amp\`ere equations:
$\det(\frac{\partial^2u}{\partial z_i\partial \overline
{z_j}})dV=\mu$, where $\mu$ was a bounded non-negative Borel
measure. Cheng and Yau in \cite{CY} solved the Dirichlet problem for
the equations $\det(\frac{\partial^2 u}{\partial z_i\partial
\overline{z_ j}})=f(z,u,\nabla u)$, $z\in\Omega$ with $f=e^u$,
$u=\infty$ on $\partial \Omega$, obtaining a solution $u\in
C^\infty(\Omega)$. In the non-degenerate case $(f>0)$, the existence
and uniqueness of a classical plurisubharmonic solution subject to
the Dirichlet boundary condition for complex Monge-Amp\`ere
equations, under suitable restrictions on $f$ and $\varphi$, have
been proved by Caffarelli, Kohn, Nirenberg and Spruck in
\cite{CKNS}.

The methods to deal with the oblique boundary value problems for
real Monge-Amp\`ere equations originated from the paper of Lions,
Trudinger and Urbas in \cite{LTU} dealing with the Neumann boundary
value problems. In 1994, Li Song-Ying in \cite{L1} considered the
Neumann problem for complex Monge-Amp\`ere equations in strictly
pseudoconvex domain in $\bf C^n$ with smooth boundary such that,
\begin{equation}
\inf_{\partial\Omega}(\gamma_0+2\lambda_2(z))>0,
\end{equation}
where $\varphi(z,u)=-\gamma_0u+\varphi(z)$, $\lambda_{2}(z)$ denotes
the smallest principal curvature of $\partial \Omega$ at $z\in
\partial\Omega$.
He obtained the existence and uniqueness of classical
plurisubharmonic solutions. The existence of generalized solutions
for
\begin{equation}\label{cma1}\det D^2u=f(x,u,Du)\ \ in \ \
\Omega\end{equation} and
\begin{equation}\label{obl}D_\beta u+\phi(x,u)=0\ \ on \ \
\partial\Omega\end{equation}
 has been established in \cite{W} by Wang Xu-jia for
general strictly oblique vector field $\beta$ under relatively weak
regularity hypotheses on the data. In \cite{U1,W} examples were
constructed showing that one cannot generally expect smoothness of
solutions of (\ref{cma}) and (\ref{ob1}), no matter how smooth
$\beta,\phi,f$ and $\partial \Omega$ are. Some structural condition
on $\beta$ is necessary. The examples in \cite{U1,W} are constructed
in such a way that $\beta$ can be made arbitrarily close to $\nu$ in
the $C^0$ norm, but not $C^1$ norm. Therefore, Urbas in \cite{U2}
proved the existence of classical solutions for real Monge-Amp\`ere
type equations subject to the semilinear oblique boundary condition
which is a $C^1$ perturbation of the Neumann boundary condition. In
1999, Li Song-Ying in \cite{L2} used different methods than Urbas in
\cite{U2} to prove the existence, uniqueness and regularity for
oblique boundary value problem to real Monge-Amp\`ere equations in a
smooth bounded strictly convex domain.

In this article we prove the existence and regularity of the
strictly plurisubharmonic solution for the oblique problem
(\ref{cma}) and (\ref{ob1}) in strictly pseudoconvex domains under
suitable conditions.

In order to prove the existence of plurisubharmonic solutions, we
have to deduce a priori estimates for the derivative of such
solutions up to second order. In the real case, Lions, Trudinger and
Urbas in \cite{LTU} used the convexity of the solution to obtain the
first order derivative bound for the oblique boundary value
problems. Unfortunately, there is no convexity for strictly
plurisubharmonic function. Inspired by the method in \cite{L1} and \cite{U2}, we
employ similar argument to obtain a priori estimates.

Now we can state our theorem as follows:
\begin{theorem}\label{thm}
Assume $\Omega$ is a bounded strictly pseudoconvex domain in $\bf
C^n$ with $\partial\Omega\in C^4$. Suppose $f\in
C^2(\overline{\Omega})$, $\varphi$ and $\beta$ satisfy the
conditions (\ref{ob2'})-(\ref{ob6}). In addition, we assume that
$\varphi\in C^{2,1}(\partial\Omega\times \bf R)$ with
\begin{equation}\label{ob3}\varphi_u(z,u)\leq -\gamma_0<0 \ on\
\partial \Omega\times \bf R,\ \ \inf_{\partial\Omega}(\gamma_0+2\lambda_2(z))>0,\end{equation} where
$\lambda_{2}(z)$ denotes the smallest principal curvature of
$\partial \Omega$ at $z\in \partial\Omega$. Then there is a positive
constant $\epsilon_0 >0$ such that
\begin{equation}\label{ob2}
\|\beta-\nu \|_{C^1(\partial \Omega)} \le \epsilon_0,
\end{equation} thus the oblique problem for
complex Monge-Amp\`ere equations (\ref{cma}) and (\ref{ob1}) has a
unique plurisubharmonic solution
 $u\in C^{2,\alpha}(\overline\Omega) $ for any $\alpha<1$.
\end{theorem}

Note that the regularity for the solution $u$ in Theorem \ref{thm}
can be improved by the linear elliptic theory in \cite{GT} if the
datas are sufficiently smooth.

The paper is organized as follows: in Section 2, we introduce some
terminology, then derive $C^0$ estimates for the solutions. In
Section 3, we study the first order derivative estimates. In section
4, we complete the second order derivative estimates. In Section 5,
we give the proof of our main theorem by using the method of
continuity.
\section{Maximum modulus estimates}

In this section, we first introduce some terminology. Then we derive
the maximum modulus estimate of the solution $u$ for the oblique
boundary value problem (\ref{cma}) and (\ref{ob1}) under some
structural conditions on $\varphi(z,u)$.

Let $z=(z_1,...,z_n)\in \bf C^n$, where $z_i=x_i+\sqrt{-1}y_i,\ (i=1,\cdots,n)$. We may write z in real coordinates
as $z=(t_1,...,t_{2n})$. Given $\xi\in \bf R^{2n}$, $D_\xi$ denotes
the directional derivative of $u$ along $\xi$. In particular,
$D_k=\frac{\partial}{\partial t_k}$. For the complex variables, we
 shall use notations: $\partial_k=\frac{\partial}{\partial z_k}$,
$\partial_{\overline k}=\frac{\partial}{\partial {\overline {z_k}}}$
and $\partial_{i \bar j}=\frac{\partial^2}{\partial z_i \partial
\overline{z_j}}$. $u^{i\bar j}$ denotes the inverse matrix of
$u_{i\bar j}$. For a compact set $X$, we let $|u|_{k,X}$ denote the
$C^k$ norm on $X$.

Let $f>0$ and $g=\log f$. Then $\log[\det(u_{i\bar j})]$ is a
concave function in $u_{i\bar j}$. We use Einstein convention. If
$u$ is a solution of (\ref{cma}), and $\xi\in \bf R^{2n}$, then
\begin{equation}\label{2.1}u^{i \bar j}\partial_{i\bar j}D_\xi u =D_\xi \log [\det
(u_{i\bar j})]=D_\xi g, \end{equation} \begin{equation}\label{2.2}
u^{i\bar j}\partial_{i\bar j}D_{\xi\xi}u\geq D_{\xi\xi}g
.\end{equation} Moreover, if we let
$\widetilde{f}=f^{\frac{1}{n}}$ and $F^{i\bar
j}=\frac{1}{n}\widetilde{f}u^{i\bar j}$, then
\begin{equation}\label{2.3}
F^{i\bar j}\partial _{i\bar j}D_\xi u= D_\xi
\widetilde{f}, F^{i\bar j}\partial_{i\bar j}D_{\xi\xi}u\geq
D_{\xi\xi}\widetilde{f},\end{equation}
 \begin{equation}\label{2.4}tr(F^{i\bar{j}})=\frac{1}{n}\widetilde{f}tr(u^{i\bar j})\geq n\frac{1}{n}\widetilde{f}\widetilde{f}^{-1}=
 1.\end{equation}

Although the argument in this
part is rather standard, for completeness, we give the maximum principle without proof:\\
\begin{lemma}\label{l1} Suppose that $L=F^{i\bar j}\partial_{i\bar j}$ is elliptic,
and $LH\geq 0 (\leq 0)\ \ in \ \Omega$ with $H\in C^2(\Omega)\cap
C^0(\overline \Omega)$. Then the maximum (minimum) of $H$ in
$\overline\Omega$ is achieved on $\partial\Omega$, that is
,\begin{equation}\sup_\Omega H=\sup _{\partial\Omega}H\ \
(\inf_\Omega H=\inf_{\partial\Omega}H).\end{equation}\end{lemma}

\begin{lemma}\label{l3} Let $\Omega$ be a bounded domain in $\bf C^n$ with $C^1$
boundary, and
 $u\in C^2(\Omega)\bigcap C^1(\overline \Omega)$ be a plurisubharmonic solution of (\ref{cma}) and (\ref{ob1}).
 Assume $\varphi$ satisfies assumption (\ref{ob4}), then
$u(z)\leq N_1 $  on $\overline\Omega$.\end{lemma}
\begin{proof} Since $u$ is plurisubharmonic, $u$ attains its maximum over
$\overline\Omega$ on $\partial\Omega$, say at $z_0\in
\partial\Omega$, then $D_\nu u(z_0)\geq0$. On the other hand, $\beta$ can be decomposed into normal and tangential part,
from (\ref{ob2'}) and the vanishing of the tangential derivative at
$z_0$, we obtain $ D_\beta u(z_0)\geq 0.$ If $u(z_0)\geq N_1,$ by
(\ref{ob4}) we have $\varphi(z_0,u(z_0))<0$, this is contradiction.
So $u(z_0)< N_1$ and the proof is complete.\end{proof}
 \begin{lemma}\label{l4} Let $u\in C^2(\Omega)\cap C^1(\overline\Omega)$ be
a plurisubharmonic solution of (\ref{cma}) and (\ref{ob1}). Assume
that $\varphi$ satisfies (\ref{ob5}), then $u(z)\geq N_2>-\infty,$
$z\in \overline \Omega$, where $N_2$ is a constant depending only on
$|f|_{0,\overline\Omega},\gamma_0,\varphi$.\end{lemma}
\begin{proof} We use the auxiliary function
\begin{equation}
h(z)=u(z)-\eta |z|^2,
\eta=|f|^{\frac{1}{n}}_{0,\overline\Omega}+1.\end{equation} We
assume that $h$ attains its minimum over $\overline\Omega$ at
$z_0\in \overline\Omega$. Now we want to show $z_0\in
\partial\Omega$.
If $z_0\in \Omega$, then
\begin{equation}
\begin{array}{ll}
0&\leq Lh(z_0)\\
&= F^{i\bar j}(z_0)\partial_{i\bar j}h(z_0)\\
&=F^{i\bar j}(z_0)\partial_{i\bar j}u(z_0)-
\eta F^{i\bar j}(z_0)\delta_{ij}\\
&=f^{\frac{1}{n}}-\eta tr(F^{i\bar j}(z_0))\\
&\leq f^{\frac{1}{n}}-\eta \\
&\leq-1\\
&<0.
\end{array}
\end{equation} This inequality leads to a contradiction, so we have shown
$z_0\in
\partial \Omega$. Without loss of generality, we can assume $u(z_0)\leq 0$,
otherwise we obtain the result. Therefore, by the oblique boundary
condition (\ref{ob1}),
\begin{equation}
\begin{array}{ll}
0 &\geq D_\beta h(z_0)\\
&=D_\beta
u(z_0)-D_\beta(\eta|z|^2)|_{z=z_0}\\
&=\varphi(z_0,u(z_0))-D_\beta(\eta|z|^2)|_{z=z_0}.
\end{array}
\end{equation}
Then
\begin{equation}
\varphi(z_0,u(z_0))\leq D_\beta(\eta|z|^2)|_{z=z_0}\leq C.
\end{equation}
Then by the condition (\ref{ob5}), we have $u(z_0)>\widetilde{N}$,
where $\widetilde{N}>-\infty$ is a constant depending on $n,\
\Omega,\ \varphi, |f|_{0,\overline{\Omega}}^{\frac{1}{n}}$. For all
$z\in\Omega$, $h(z)\geq h(z_0)$, which leads to,
\begin{equation}
\begin{array}{ll}
u(z)
&=h(z)+\eta|z|^2\\
&\geq h(z_0)+\eta|z|^2\\
&=u(z_0)+\eta|z|^2-\eta|z_0|^2\\
&\geq \widetilde{N}+\eta|z|^2-\eta|z_0|^2\\
&\geq N_2\\
&>-\infty,
\end{array}
\end{equation}
where $N_2$ depends on $n,\ \Omega,\ \varphi,
|f|_{0,\overline{\Omega}}^{\frac{1}{n}}$. The proof of Lemma
\ref{l4} is complete.\end{proof}
\begin{theorem}\label{thm2} Let $\Omega$ be a bounded strictly pseudoconvex domain
with $C^1$ boundary. Assume that $f$ is non-negative, and $\varphi
$ satisfies
(\ref{ob4})-(\ref{ob5}). If $u\in C^2(\Omega)\cap
C^1(\overline\Omega)$
 is a plurisubharmonic solution of (\ref{cma}) and (\ref{ob1}),  \ then
 $|u|_{0,\overline \Omega}\leq C$, where $C$ is a constant depending only on
 $n,\Omega,|f|_{0,\overline \Omega},
 \gamma_0,\varphi$.\end{theorem}
\begin{remark} In the proof of Lemma \ref{l4}, the estimate on
$|u|_{0,\overline\Omega}$ is independent of the lower bound of $f$.
So when $f\geq 0$, we can obtain the same estimate by considering
$f_\epsilon=f+\epsilon, \epsilon>0$ first, then let
$\epsilon\rightarrow 0^+$ to complete the proof of theorem
\ref{thm2}. If a similar argument is needed in the following
sections, we will not repeat again.\end{remark}
\begin{remark}
We remark that $\widetilde{C_l}$ and $C_l$ $(l=1,2,\cdots)$denote the constants depending on the known datas. As usual,
 constants may change from line to line in the context.
\end{remark}
\section{Gradient estimates}

In this section, we follow the idea in \cite{L1} and \cite{U2} to derive gradient
estimates.
\begin{theorem}\label{ge}
Let $\Omega$ be a bounded strictly pseudoconvex domain in $\bf C^n$
with $C^3$ boundary. Assume $\beta$, $\varphi\in C^{1,1}(\partial
\Omega\times \bf R)$ and $f\in C^1(\overline{\Omega})$ satisfy
(\ref{ob2'})-(\ref{ob6}). In addition $\varphi$ satisfies
(\ref{ob3}). Then
\begin{equation}\label{3.1'}|u|_{1,\overline{\Omega}}\leq
C,\end{equation} where $C$ is a constant depending only on $\
|f^{\frac{1}{n}}|_{1,\overline\Omega},\ \gamma_0,\ \Omega,\
\beta\ and \ |\varphi|_{C^{1,1}(\partial\Omega)}$.
\end{theorem}
\begin{proof}
In order to prove $|u|_{1,\overline\Omega}\leq C$, from Theorem
\ref{thm2} it suffices to prove \begin{equation}\label{g1} D_\xi
u(z)\leq C,\ \ \forall\xi\in S^{2n-1},
\end{equation}
where $S^{2n-1}$ is the unit sphere in $R^{2n}$. We still use
$\varphi$ to denote a $C^1$ extension of $\varphi$ from
$\partial\Omega$ to $\overline\Omega.$

First, we reduce the gradient estimates to the boundary by choosing the auxiliary function $R(z,\xi)=D_{\xi}u(z)+\eta_1|z|^2$, where $\eta_1=|f^{\frac{1}{n}}|_{1,\overline{\Omega}}+1$. By calculation, we have $LR>0$. By the maximum principle,
\begin{equation}\label{3.3}
\begin{array}{ll}
\max_{\overline{\Omega}}D_\xi u(z)&\leq \max_{\overline{\Omega}}(D_\xi u(z)+\eta_1|z|^2)\\
&\leq\max_{\partial\Omega}(D_\xi u(z)+\eta_1|z|^2)\\
&\leq\max_{\partial\Omega}D_\xi u(z)+\eta_1 diam(\Omega)^2.
\end{array}
\end{equation}

So the rest task is the estimation of $D_\xi u(z)$ on $\partial\Omega$.
From the condition (\ref{ob1}), when $z\in\partial\Omega$,
 \begin{equation}\label{go}
|D_{\beta} u(z)|=|\varphi(z,u(z))|\leq C_1.
 \end{equation}

At any boundary point, any direction $\xi$ can be written in terms
of a tangential component $\tau(\xi)$
 and a component in the direction $\beta$, namely
 \begin{equation}
 \xi=\tau(\xi)+\frac{(\nu\cdot\xi)}{(\beta\cdot\nu)}\beta,
 \end{equation}
 where $\tau(\xi)=\xi-(\nu\cdot\xi)\nu - \frac{(\nu\cdot\xi)}{(\beta\cdot\nu)} \beta^T$ and $\beta^T=\beta-(\beta\cdot\nu)\nu$.
We compute under the condition (\ref{ob2'}),
\begin{equation}\label{3.6}
\begin{array}{ll}
|\tau(\xi)|^2&=1-(1-\frac{|\beta^T|^2}{|(\beta\cdot\nu)|^2})(\nu\cdot\xi)^2-2\frac{(\nu\cdot\xi)(\beta^T\cdot\xi)}{(\beta\cdot\nu)}\\
&\leq 1+\frac{2}{\beta_0}|\beta^T|.
\end{array}
\end{equation}
Thus, if we get the tangential derivative estimate of $u$ on
$\partial\Omega$, from (\ref{go}), for any $\xi\in S^{2n-1}$, we
have
\begin{equation}\label{3.7}
\begin{array}{ll}
D_{\xi}u(z)&=D_{\tau(\xi)}u(z)+\frac{(\nu\cdot\xi)}{(\beta\cdot\nu)}D_\beta u(z)\\
&\leq D_{\tau(\xi)}u(z)+C_2\\
&=|\tau(\xi)|D_{\frac{\tau(\xi)}{|\tau(\xi)|}}u(z)+C_2\\
&\leq \sqrt{(1+\frac{2}{\beta_0}|\beta^T|)}C_3+C_2\\
&\leq C_4.
\end{array}
\end{equation}

Next the key task is to get the bound for the tangential derivative
of $u$ on $\partial \Omega$. We assume that the maximum tangential first order derivative on $\partial\Omega$ is attained
at a boundary point which may take to be the origin, in a tangential direction which may take to be $e_1=(1,\underbrace{0,\cdots,0}_{2n-1})$, $x_n$ is the inner normal vector at 0. Thus
\begin{equation}\label{3.8}
D_1u(0)=\sup_{z\in \partial\Omega,\tau is\ unit\ tangential\  at z }D_\tau u(z).
\end{equation}
And without loss of generality, we can assume $D_1u(0)>0$. Otherwise, we finish the proof.

We choose the auxiliary function
\begin{equation}
Q(z)=\frac{D_{1}u(z)}{D_1u(0)}+G|z|^2-Ax_n-1,
\end{equation} and the domain $S_\mu=\{z\in\Omega:x_n\leq \mu\}$, where $G,A $ and $\mu$ are constants to be fixed.

From (\ref{2.3}) we have $F^{i\bar j}\partial_{i\bar
j}\frac{D_1 u(z)}{D_1 u(0)}=\frac{D_1 f^{\frac{1}{n}}}{D_1u(0)}\geq -\frac{C_5}{D_1u(0)},$
and
\begin{equation}
F^{i\bar j}\partial_{i\bar j}G|z|^2=G tr(F^{i\bar j})\geq G.\\
\end{equation}
Thus, if we take
\begin{equation}\label{c1}G\geq\frac{C_5}{D_1u(0)},\end{equation} we have
\begin{equation}
LQ= F^{i\bar j}Q_{i\bar j}\geq0\ \ \ in\  S_\mu.
\end{equation}

Then we consider the estimates of $Q(z)$ on the boundary of $S\mu$. On $\partial\Omega\cap\overline{S_\mu}$ near $0$,
there is a constant $a>0$ such that(see  in \cite{CKNS}, Lemma 1.3)
\begin{equation}
x_n\geq a|z|^2.
\end{equation}
Thus \begin{equation}
\begin{array}{ll}
Q(z)&\leq 1+\frac{C_6|z|^2}{D_1u(0)}+G|z|^2-Ax_n-1\\
&\leq (G+\frac{C_6}{D_1u(0)}-Aa)|z|^2.
\end{array}
\end{equation}
If we take \begin{equation}\label{c2}Aa\geq G+\frac{C_6}{D_1u(0)},\end{equation} then $Q(z)\leq 0$ on $\partial\Omega\cap\overline{S_\mu}$ near $0$.

On the other hand, from (\ref{3.3}),(\ref{3.7}) and (\ref{3.8}) we have
\begin{equation}
\begin{array}{ll}
D_1u(z)&\leq \max_{\partial\Omega}D_1 u(z)+\eta_1 diam(\Omega)^2\\
&\leq\max_{\partial\Omega}\sqrt{(1+\frac{2}{\beta_0}|\beta^T|)}D_1u(0)+C_7.
\end{array}
\end{equation}
Then on $\{x_n=\mu\}\cap \overline{S_\mu}$,
\begin{equation}
\begin{array}{ll}
Q(z)&\leq \sqrt{(1+\frac{2}{\beta_0}|\beta^T|)}+\frac{C_7}{D_1u(0)}+G|z|^2-A\mu-1\\
&\leq \sqrt{(1+\frac{2}{\beta_0}|\beta^T|)}+\frac{C_7}{D_1u(0)}+G|z|^2-A\mu-1\\
&\leq (1+\frac{2}{\beta_0}|\beta^T|)+\frac{C_7}{D_1u(0)}+G|z|^2-A\mu-1\\
&\leq \frac{2}{\beta_0}|\beta^T|+\frac{C_7}{D_1u(0)}+GC\mu-A\mu.
\end{array}
\end{equation}
If \begin{equation}\label{c3}\frac{2}{\beta_0}|\beta^T|+GC\mu+\frac{C_7}{D_1u(0)}\leq A\mu,\end{equation} then $Q\leq 0$ on $\{x_n=\mu\}\cap \overline{S_\mu}$.

We now proceed to fix $G$ and $\mu$, depending on $A$ which will be fixed later. First we fix $G>0$ so small that
\begin{equation}
CG\leq \frac{1}{2}A\ and \ G\leq\frac{1}{2}Aa,
\end{equation}
then fix $\mu\in(0,1)$. Then (\ref{c1}), (\ref{c2}) and (\ref{c3}) will be hold whenever
\begin{equation}\label{gzc}
\begin{array}{ll}
&\frac{C_5}{D_1u(0)}\leq G,\\
& \frac{C_6}{D_1u(0)}\leq\frac{Aa}{2},\\
& \frac{2}{\beta_0}|\beta^T|+\frac{C_7}{D_1u(0)}\leq \frac{1}{2}A\mu.
\end{array}
\end{equation}
By the maximum principle we have $Q\leq 0$ in $S_\mu$, and since $Q(0)=0$, then
\begin{equation}D_\beta Q(0)\geq 0.\end{equation}

From (\ref{ob3})
\begin{equation}
\begin{array}{ll}
 &D_\beta D_1 u(0)\\
  &=D_1D_\beta u(0)-\sum_{k=1}^{2n}(D_1\beta_k)D_ku(0)\\
  &=D_1\varphi(z,u)(0)-\sum_{k=1}^{2n}(D_1\beta_k)D_ku(0)\\
  &=\varphi_1(0,u(0))+\varphi_u(0,u(0))D_1u(0)-\sum_{k=1}^{2n}(D_1\beta_k)D_ku(0)\\
  &\leq -\gamma_0D_1u(0)-\sum_{k=1}^{2n}(D_1\beta_k)D_ku(0)+C,
  \end{array}
\end{equation}
where $\beta_k$ is the component of $\beta$.
   Because
  $\Omega$ is a bounded strictly pseudoconvex domain, there is a
  strongly plurisubharmonic defining function $r(z)$ for $\Omega$
  satisfying $|\nabla r|=1$ on $\partial\Omega=\{z\in {\bf
  C^n}:r(z)=0\}$. And there is an orthogonal matrix $U$ over $R^{2n-1}$ such that
\begin{equation}\label{U}
U^t [\frac{\partial^2r(0)}{\partial t_k \partial
  t_l}]_{(2n-1)\times(2n-1)}
  U=diag(\lambda_2(0),\cdots,\lambda_{2n}(0)),\end{equation}where
  $\lambda_{2n}(0)\geq\cdots\geq\lambda_2(0)$. Since $\Omega$ is a bounded domain with $C^2$ boundary, there are constants $\Lambda$ and $\lambda$
  such that $\Lambda>\sup\{\lambda_{2n}(z):z\in\partial\Omega\}$
  and $\lambda<\inf\{\lambda_2(z):z\in\partial\Omega\}$.

To handle the $\sum_{k=1}^{2n}(D_1\beta_k)D_ku(0)$, we express $\nu$ in terms of
$\beta$ and tangential components,
\begin{equation}
D_{\nu}u(0)=\sum_{k=1}^{2n-1}(-\frac{\beta_k}{\beta_{2n}})D_ku(0)+\frac{1}{\beta_{2n}}D_\beta
u(0),
\end{equation}
and $\nu_k(z)=\frac{\partial r}{\partial t_k}(z)$, $z\in\partial\Omega$, $\nu_k$ is the component of $\nu$. Then
\begin{equation}
\begin{array}{ll}
\sum_{k=1}^{2n}(D_1\beta_k)D_ku(0)&=\sum_{k=1}^{2n}(D_1\beta_k-D_1\nu_k)D_k
u(0)+\sum_{k=1}^{2n}(D_1\nu_k)D_ku(0)\\
&=\sum_{k=1}^{2n-1}\frac{\partial^2 r}{\partial t_k\partial
  t_1}D_ku(0)+\sum_{k=1}^{2n-1}D_1(\beta_k-\nu_k)D_k u(0)\\
  &+[\frac{\partial^2 r}{\partial t_{2n}\partial
  t_1}+D_1(\beta_{2n}-\nu_{2n})]D_\nu u(0)\\
  &=\sum_{k=1}^{2n-1}\frac{\partial^2 r}{\partial t_k\partial
  t_1}D_ku(0)+\sum_{k=1}^{2n-1}D_1(\beta_k-\nu_k)D_k u(0)\\
  &+[\frac{\partial^2 r}{\partial t_{2n}\partial
  t_1}+D_1(\beta_{2n}-\nu_{2n})][\sum_{k=1}^{2n-1}(-\frac{\beta_k}{\beta_{2n}})D_ku(0)+\frac{1}{\beta_{2n}}D_\beta
u(0)].
\end{array}
\end{equation}
Thus
\begin{equation}\label{b}
  \begin{array}{ll}
  & D_\beta D_1u(0)\\
  &\leq C-\gamma_0D_{1}u(0)-\sum_{k=1}^{2n-1}\frac{\partial^2 r}{\partial t_k\partial
  t_1}D_ku(0)-\sum_{k=1}^{2n-1}D_1(\beta_k-\nu_k)D_ku(0)\\
  &-[\frac{\partial^2 r}{\partial t_{2n}\partial
  t_1}+D_1(\beta_{2n}-\nu_{2n})][\sum_{k=1}^{2n-1}(-\frac{\beta_k}{\beta_{2n}})D_{k}u(0)+\frac{1}{\beta_{2n}}D_{\beta}u(0)]\\
  &\leq C-\gamma_0D_{1}u(0)-[\sum_{k=1}^{2n-1}\frac{\partial^2 r}{\partial t_k\partial
  t_1}+\sum_{k=1}^{2n-1}D_1(\beta_k-\nu_k)\\
  &+\sum_{k=1}^{2n-1}(\frac{\partial^2 r}{\partial t_{2n}\partial
  t_1}+D_1(\beta_{2n}-\nu_{2n}))(-\frac{\beta_k}{\beta_{2n}})]D_ku(0)\\
  &=C-\gamma_0D_{1}u(0)-3\Lambda
  D_{1}u(0)+\sum_{k=1}^{2n-1}[\Lambda\delta_{k1}-\frac{\partial^2 r}{\partial t_k\partial
  t_1}]D_ku(0)\\
  &+\sum_{k=1}^{2n-1}[\Lambda\delta_{k1}-D_1(\beta_k-\nu_k)]D_ku(0)\\
  &+\sum_{k=1}^{2n-1}[\Lambda\delta_{k1}+
  (\frac{\partial^2 r}{\partial t_{2n}\partial
  t_1}+D_1(\beta_{2n}-\nu_{2n}))(\frac{\beta_k}{\beta_{2n}})]D_ku(0).
  \end{array}
  \end{equation}
Let \begin{equation}
A_1=\sum_{k=1}^{2n-1}[\Lambda\delta_{k1}-\frac{\partial^2
r}{\partial t_k\partial
  t_1}]D_ku(0),
\end{equation}
\begin{equation}
A_2=\sum_{k=1}^{2n-1}[\Lambda\delta_{k1}-D_1(\beta_k-\nu_k)]D_ku(0),
\end{equation}
\begin{equation}
A_3=\sum_{k=1}^{2n-1}[\Lambda\delta_{k1}+
  (\frac{\partial^2 r}{\partial t_{2n}\partial
  t_1}+D_1(\beta_{2n}-\nu_{2n}))(\frac{\beta_k}{\beta_{2n}})]D_ku(0).
\end{equation}
By using an argument similar to that given in the
  proof of Theorem 3.1 in \cite{L1}, we want to obtain the
  relationship between $A_1,\ A_2 $ and $A_3$ with $D_{1}u(0)$.

First, by the conclusion in \cite{L1} directly, we set
\begin{equation}
\hat{\tau_1}=(e'_1[\Lambda I_{2n-1}-(\frac{\partial^2 r(0)}{\partial
t_k\partial t_l})^{2n-1}_{k,l=1}],0),
\end{equation}
where $e_1'=(1,\underbrace{0,\cdots,0}_{2n-2})$. Since the matrix
$[\Lambda I_{2n-1}-(\frac{\partial^2 r(0)}{\partial t_k\partial
t_l})^{2n-1}_{k,l=1}]$ is non-negative definite with maximum
eigenvalue $\Lambda-\lambda_2(0)$, we have $|\hat{\tau_1}|\leq
\Lambda-\lambda_2(0)$. Thus by (\ref{3.8}) and $\hat{\tau_1}$ is the tangential direction at 0,
\begin{equation}\label{b1}
\begin{array}{ll}
A_1&=D_{\hat{\tau_1}}u(0)\\
&=|\hat{\tau_1}|D_{\frac{\hat{\tau_1}}{|\hat{\tau_1}|}}u(0)\\
&\leq (\Lambda-\lambda_2(0))D_1u(0).
\end{array}
\end{equation}
Then, set
\begin{equation}
M_{2n-1}=\left(
      \begin{array}{cccc}
        D_1(\beta_1-\nu_1) & D_1(\beta_2-\nu_2)&\cdots & D_1(\beta_{2n-1}-\nu_{2n-1})\\
        D_1(\beta_2-\nu_2) & 0 &\cdots & 0 \\
        \vdots &\vdots &\ddots&0\\
        D_1(\beta_{2n-1}-\nu_{2n-1})& 0&\cdots & 0\\
      \end{array}
    \right),
\end{equation}
\begin{equation} \hat{\tau_2}=(e'_1[\Lambda I_{2n-1}-M_{2n-1}],0).
\end{equation}
 Thus
\begin{equation}
A_2=D_{\hat{\tau_2}}
u(0)=|\hat{\tau_2}|D_{\frac{\hat{\tau_2}}{|\hat{\tau_2}|}}u(0).
\end{equation}
By calculation, the eigenvalues of matrix $M_{2n-1}$ are
\begin{equation}
\pi_1=0,\ \
\pi_{2,3}=\frac{D_1(\beta_1-\nu_1)\pm\sqrt{(D_1(\beta_1-\nu_1))^2+4[\sum_{i=2}^{2n-1}(D_1(\beta_i-\nu_i))^2]}}{2},
\end{equation}
where $\pi_1$ is (2n-3) multiple eigenvalues. In the condition
(\ref{ob2}), we can take
$\epsilon_0\leq\frac{\Lambda}{1+\sqrt{8n-7}}$ so that the matrix
$\Lambda I_{2n-1}-M_{2n-1}$ is non-negative definite with
$|\hat{\tau_2}|\leq (\Lambda+\frac{1+\sqrt{8n-7}}{2}\epsilon_0)$.
Thus
\begin{equation}\label{b2'}
A_2\leq (\Lambda+\frac{1+\sqrt{8n-7}}{2}\epsilon_0)D_{1}u(0).
\end{equation}
At last, we consider $A_3$.
\begin{equation}
A_3=\sum_{k=1}^{2n-1}[\Lambda\delta_{k1}+
  (\frac{\partial^2 r}{\partial t_{2n}\partial
  t_1}+D_1(\beta_{2n}-\nu_{2n}))(\frac{\beta_k}{\beta_{2n}})]D_ku(0).
\end{equation}
Set
\begin{equation}N(0)=[\frac{\partial^2r}{\partial t_{2n}\partial
  t_1}+D_1(\beta_{2n}-\nu_{2n})](0),\end{equation}
here $N(0)$ has uniform upper bound, i.e. $|N(0)|\leq N$, where $N$ is a constant independent of $0$. Set
  \begin{equation}
 G_{2n-1}=\left(
                  \begin{array}{cccc}
                  \beta_1 & \beta_2 & \cdots & \beta_{2n-1} \\
                    \beta_2& 0 & \cdots & 0 \\
                    \vdots & \vdots & \ddots &\vdots\\
                    \beta_{2n-1} & 0 &\cdots & 0\\
                  \end{array}
                \right),
  \end{equation}
  \begin{equation}
  \hat{\tau_3}=(e_1'[\Lambda I_{2n-1}+\frac{N(z_0)}{\beta_{2n}}G_{2n-1}],0).
  \end{equation}
  So, if we take $\epsilon_0\leq\frac{\Lambda|\beta_{2n}|}{N(1+\sqrt{8n-7})}$, we have
  \begin{equation}\label{b3'}
  \begin{array}{ll}
  A_3=D_{\hat{\tau_3}} u(0) =|\hat{\tau_3}|D_{\frac{\hat{\tau_3}}{|\hat{\tau_3}|}}u(0)\leq
 (\Lambda+|\frac{N}{\beta_{2n}}|\frac{1+\sqrt{8n-7}}{2}\epsilon_0)D_{1}u(0).
  \end{array}
  \end{equation}
Thus by inserting (\ref{b1}), (\ref{b2'}) and (\ref{b3'}) into
(\ref{b}) we get
\begin{equation}\label{3.37}
\begin{array}{ll}
&D_\beta D_1u(0)\\
&\leq C-\gamma_0D_{1}u(0) -3\Lambda
  D_{1}u(0)+(\Lambda-\lambda_2(0))D_{1}u(0)\\
  &+(\Lambda+\frac{1+\sqrt{8n-7}}{2}\epsilon_0)D_{1}u(0)
+(\Lambda+|\frac{N}{\beta_{2n}}|\frac{1+\sqrt{8n-7}}{2}\epsilon_0)D_{1}u(0)\\
&=C-[\gamma_0+\lambda_2(0)-\frac{(1+\sqrt{8n-7})\epsilon_0}{2}-\frac{(1+\sqrt{8n-7})|\frac{N}{\beta_{2n}}|\epsilon_0}{2}]D_{1}u(0).
  \end{array}
\end{equation}
Again, we take $\epsilon_0\leq
\frac{(\gamma_0+\lambda_2(0))}{[(1+\sqrt{8n-7})(1+|\frac{N}{\beta_{2n}}|)]}$,
then
\begin{equation}
{\gamma_0+\lambda_2(0)}\geq
(1+\sqrt{8n-7})\epsilon_0+(1+\sqrt{8n-7})|\frac{N}{\beta_{2n}}|\epsilon_0.
\end{equation}
Ultimately, by (\ref{ob2'}),
$\beta_{2n}(0)=\beta(0)\cdot\nu(0)\geq \beta_0$, we choose
\begin{equation}
\begin{array}{ll}
 \epsilon_0&\leq
\min\{\frac{\Lambda}{1+\sqrt{8n-7}},\frac{\Lambda\beta_0}{N(1+\sqrt{8n-7})},\frac{(\gamma_0+\lambda_2(z_0))}{[(1+\sqrt{8n-7})(1+\frac{N}{\beta_0})]}\}\\
&\leq\min\{\frac{\Lambda}{1+\sqrt{8n-7}},\frac{\Lambda|\beta_{2n}|}{N(1+\sqrt{8n-7})},
\frac{(\gamma_0+\lambda_2(z_0))}{[(1+\sqrt{8n-7})(1+|\frac{N}{\beta_{2n}}|)]}\},
\end{array}
\end{equation}
then
\begin{equation}\label{3.40} [\gamma_0+\lambda_2(0)-\frac{(1+\sqrt{8n-7})\epsilon_0}{2}-\frac{(1+\sqrt{8n-7})|\frac{N}{\beta_{2n}}|\epsilon_0}{2}]\geq \frac{\sigma_0}{2},
\end{equation}
where $\sigma_0=\inf_{\partial\Omega}(\gamma_0+\lambda_2)$.
Then by inserting (\ref{3.40}) into (\ref{3.37}),
\begin{equation}
\begin{array}{ll}
0&\leq D_\beta Q(0)\\
&\leq \frac{C-\frac{\sigma_0}{2}D_1u(0)}{D_1u(0)}-AD_\beta x_n(0).
\end{array}
\end{equation}
\begin{equation}
[\frac{\sigma_0}{2}-A\beta\cdot\nu(0)]D_1u(0)\leq C.
\end{equation}
Finally we fix $A$ so small that
\begin{equation}
A\beta\cdot\nu\leq\frac{\sigma_0}{4}, on\ \partial\Omega,
\end{equation}
 we have for any tangential vector  $D_1
u(0)\leq \frac{\sigma_0}{4}.$ We finish the proof of Theorem \ref{ge}.

\end{proof}

\section{Second order derivative estimates }
In this section, we aim to derive the second order derivative
estimates. First, we reduce the interior second order derivatives to
the boundary. Then, we derive the second order derivative estimates
on the boundary. Finally, from these estimates we have Theorem
\ref{thm4}.

\begin{theorem}\label{thm3} Let $\Omega$ be a bounded strictly pseudoconvex domain in
$\bf C^n$ with $C^4$ boundary. Assume $\beta$, $\varphi\in
C^{2,1}(\partial \Omega\times \bf R)$ and $f\in
C^2(\overline{\Omega})$ satisfy (\ref{ob2'})-(\ref{ob6}). In
addition $\varphi$ satisfies (\ref{ob3}). Then
\begin{equation}\label{3.1}|D^2 u|_{0,\overline{\Omega}}\leq
\widetilde{C},\end{equation} where $\widetilde{C}$ is a constant depending only on
$|u|_{1,\overline\Omega},\
|f^{\frac{1}{n}}|_{C^{2}(\overline\Omega)},\ \gamma_0,\ \Omega,
\beta\ and \ |\varphi|_{C^{2,1}(\partial\Omega)}$.\end{theorem}
First, we reduce the second order derivative estimates to the boundary by choosing the auxiliary function $R_1(z,\xi)=D_{\xi\xi}u(z)+\eta_2|z|^2$, where $\eta_2=|f^{\frac{1}{n}}|_{C^{2}(\overline{\Omega})}+1$.

 We have
 \begin{equation}F^{i\bar j}\partial_{i\bar j}D_{\xi\xi}u
 \geq D_{\xi\xi}\widetilde{f}\geq -|f^{\frac{1}{n}}|_{C^{2}(\overline\Omega)},\end{equation}
 thus $LR_1=F^{i\bar j}\partial_{i\bar j}R_1>0$. By the maximum principle,
\begin{equation}\label{4.3}
\begin{array}{ll}
\max_{\overline{\Omega}}D_{\xi\xi} u(z)&\leq \max_{\overline{\Omega}}(D_{\xi\xi} u(z)+\eta_2|z|^2)\\
&\leq\max_{\partial\Omega}(D_{\xi\xi} u(z)+\eta_2|z|^2)\\
&\leq\max_{\partial\Omega}D_{\xi\xi} u(z)+\eta_2 diam(\Omega)^2.
\end{array}
\end{equation}

Next, we give some lemmas below.
\begin{lemma}\label{l31} We reformulate (\ref{3.1}) as follow,
\begin{equation}\label{3.2}D_{\xi\xi}u(z)\leq C,\ \ \ \xi\in \bf R^{2n},\end{equation}whenever $u$ is
subharmonic.\end{lemma}

\begin{lemma}\label{l32} $|D_{\beta\tau }u|\leq C$ on $\partial \Omega$  for any unit tangential vector field
$\tau$.\end{lemma}
\begin{proof}
By applying the tangential gradient operator to the boundary
condition we obtain $|D_{\beta\tau }u|\leq C$ on $\partial \Omega$
for any tangential vector $\tau$.
\end{proof}
\begin{lemma}\label{l33}$| D_{\beta\nu}u|\leq C$ on $\partial
\Omega$.\end{lemma}
\begin{proof}
Without loss of generality, we set $0\in
\partial \Omega$, $x_n$ is inner normal vector at 0.
Case 1: Suppose that $Du(0)=0$.

Now we use the auxiliary function \begin{equation}W=\pm D_\beta u\mp
\varphi(z,u)+|Du|^2+\widetilde{K}|z|^2-\widetilde{G}x_n.\end{equation}
Now we compute $LW$.
\begin{equation}
\begin{array}{ll}
LD_\beta u &=F^{i\bar j}\partial _i \partial _{\bar
j}\beta_k D_k u+F^{i\bar j}(\partial _{\bar j}\beta_k)(\partial _i
D_k u)+ F^{i\bar j}(\partial _i \beta_k)(\partial _{\bar j}D_k u)+
F^{i\bar j}\beta_k\partial_i\partial_{\bar j}D_k u\\
&=F^{i\bar j}\partial _i \partial _{\bar j}\beta_k D_k u+F^{i\bar
j}(\partial _{\bar j}\beta_k)(\partial _i D_k u)+ F^{i\bar
j}(\partial _i \beta_k)(\partial _{\bar j}D_k
u)+\beta_kD_k\widetilde{f}.
\end{array}
\end{equation}
From the Cauchy inequality and the positivity of the matric
$(F^{i\bar j})$, we have \begin{equation}|F^{i\bar j}(\partial
_{\bar j}\beta_k)(\partial _i D_k u)|\leq(F^{i\bar
j}(\partial_i\beta_k\partial _{\bar
j}\beta_k))^{\frac{1}{2}}(F^{i\bar j}(\partial _i D_k
u\partial_{\bar j}D_k u))^{\frac{1}{2}},\end{equation}
\begin{equation}|F^{i\bar j}(\partial _i\beta_k)(\partial _{\bar j} D_k
u)|\leq(F^{i\bar j}(\partial_i\beta_k\partial _{\bar
j}\beta_k))^{\frac{1}{2}}(F^{i\bar j}(\partial _i D_k
u\partial_{\bar j}D_k u))^{\frac{1}{2}}. \end{equation}Hence
\begin{equation}|LD_\beta u|\leq \widetilde{C_1}tr (F^{i\bar
j})+F^{i\bar j}(\partial_i D_k u)(\partial _{\bar j}D_k
u)+C,\end{equation} where $\widetilde{C_1}$ depends on
$|u|_{1,\overline\Omega}$, $\beta$ and $f^{\frac{1}{n}}$. We also
have
\begin{equation}
\begin{array}{ll}
L\varphi(z,u)&=F^{i\bar j}\partial _i \partial _{\bar
j}\varphi(z,u)\\
&=F^{i\bar j}[\varphi_{i\bar j}+\varphi_{u\bar
j}u_i+\varphi_{ui}u_{\bar j}+\varphi_{uu}u_iu_{\bar
j}+\varphi_uu_{i\bar j}]\\
&=F^{i\bar j}[\varphi_{i\bar j}+\varphi_{u\bar
j}u_i+\varphi_{ui}u_{\bar j}+\varphi_{uu}u_iu_{\bar
j}]+\widetilde{f}\varphi_u\\
&\leq \widetilde{C_2}tr(F^{i\bar j})+C,
\end{array}
\end{equation} where $\widetilde{C_2}$ depends on $\varphi$, $f^{\frac{1}{n}}$ and $|u|_{1,\overline\Omega}$.
\begin{equation}
\begin{array}{ll}
L|Du|^2 &=2F^{i\bar j}(\partial _i D_k u\partial_{\bar
j}D_k u+D_k u\partial_{i\bar
j}D_ku)\\
&=2F^{i\bar j}(\partial _iD_k u)(\partial_{\bar j}D_k u)+2D_k uD_k
\widetilde{f}\\
&\geq 2F^{i\bar j}(\partial _iD_k u)(\partial_{\bar j}D_k u)-C,
\end{array}
\end{equation}
where $C$ depends on $|u|_{1,\overline\Omega}$ and
$f^{\frac{1}{n}}$.
\begin{equation}L|z|^2=tr(F^{i\bar j}).\end{equation} Therefore, we obtain
\begin{equation}LW\geq(\widetilde{K}-\widetilde{C_1}-\widetilde{C_2})tr(F^{i\bar j})-\widetilde{C_3}\geq\widetilde{K}-\widetilde{C_4},\end{equation} we
can choose $\widetilde{K}>\widetilde{C_4}$ such that $LW\geq 0.$

Finally, we use a standard barrier argument for $D_\nu W$ on
$\partial\Omega$(see in \cite{CKNS}, Lemma 1.3 or \cite{GT},
Corollary 14.5). Let $S_{\mu_1}=\{z\in \Omega|x_n \leq \mu_1\}. $ On $\partial S_{\mu_1}\cap\Omega $ if
$\widetilde{G}$ is sufficiently large, we get $W\leq0$. On $\partial
S_{\mu_1}\cap\partial\Omega$, we have for some $a>0$, $x_n\geq
a|z|^2$. Then we can choose $\widetilde{G}$ large enough such that $W\leq 0 $ on
$\partial S_{\mu_1}\cap\partial\Omega$. Since $W(0)=0$, by the
maximum principle, we obtain $W_{x_n}(0)\leq 0$, then $|D_{\beta
x_n}u|\leq \widetilde{G}+|\varphi|_{C^{1,1}(\partial \Omega)}.$

Case 2: If $Du(0)\neq0$, we can take \begin{equation}W=\pm D_\beta u\mp
\varphi(z,u)+|Du(z)-Du(0)|^2+\widetilde{K}|z|^2-\widetilde{G}x_n,\end{equation} the proof
above is still valid.
\end{proof}
\begin{remark}Combining Lemma \ref{l32}  with Lemma \ref{l33} , for
any direction $\xi$, we obtain $|D_{\beta\xi}u|\leq C$\ on $
\partial \Omega$. In particular $|D_{\beta\beta}u|\leq C$ on
$\partial \Omega.$
\end{remark}

Since we have the bounds for $|D_{\tau\beta}u|$ and
$|D_{\beta\beta}u|$ on $\partial\Omega$, in order to finish the proof of Theorem
\ref{thm3}, the remaining task is to get the bounds for the
tangential second derivatives of $u$ on $\partial\Omega$. We now assume that the maximum tangential second order derivative on
$\partial\Omega$ is attained at a boundary point which may take to be the origin, in tangential direction which we may take to be $e_1$, and $x_n$ is the inner
normal vector at 0. Thus
\begin{equation}\label{4.15}
D_{11}u(0)=\sup_{z\in\partial\Omega,\tau\  is\ unit\  tangential\  at\  z}D_{\tau\tau}u(z).
\end{equation}
Without loss of generality, we can take $D_{11}u(0)>0$. Otherwise, we finish the proof.

As for any direction $\xi$, on $\partial\Omega$,
\begin{equation}\label{4.16}
\begin{array}{ll}
D_{\xi\xi}u&=D_{\tau(\xi)\tau(\xi)}u+2\frac{(\nu\cdot\xi)}{(\beta\cdot\nu)}D_{\tau(\xi)\beta}u
+\frac{(\nu\cdot\xi)^2}{(\beta\cdot\nu)^2}D_{\beta\beta}u, \\
&\leq (1+\frac{2}{\beta_0})D_{11}u(0)+C\ \ on \
\
\partial\Omega.
\end{array}
\end{equation}

We introduce the tangential gradient operator
$\delta=(\delta_1,\cdots, \delta_{2n-1})$, where
$\delta_i=(\delta_{ij}-\nu_i\nu_j)D_j$. Applying this tangential
operator to the boundary condition (\ref{ob1}), we have
\begin{equation}
D_ku\delta_i \beta_k + \beta_k \delta_i D_ku=\delta_i\varphi,\quad
{\rm on} \  \partial \Omega,
\end{equation}
then
\begin{equation}
D_{\tau\beta}u=-D_ku(\delta_i\beta_k)\tau_i+\delta_i\varphi\tau_i.
\end{equation}

For the case $\xi=e_1$ in (\ref{4.16}), we have
\begin{equation}\label{4.19}
D_{11}u(z)\leq (1-\frac{2\nu_1\beta_1^T}{\beta\cdot\nu})D_{11}u(0)+\frac{2\nu_1}{\beta\cdot\nu}D_{\tau(e_1)\beta}+\frac{\nu^2_1}{(\beta\cdot\nu)^2}D_{\beta\beta}u(z).
\end{equation}

 Similarly to the real case in \cite{U2}, let $S_{\mu_2}=\{z\in\Omega:x_n\leq \mu_2\}$. We construct a function
 $$H(z)=\frac{D_{11}u(z)-V(z)}{D_{11}u(0)}+\frac{2\nu_1\beta_1^T}{\beta\cdot\nu}+\widehat{B}|Du(z)|^2+\widehat{G}|z|^2-\widehat{A}x_n-1,$$
 $\widehat{G}$, $\widehat{A}$, $\widehat{B}$ and $\mu_2$ are constants to be fixed.
$V(z)$ is a linear function with respect to $Du$, such that
 \begin{equation}
 V(z)=a_k(z)D_ku+b(z),\ \ in\ \  \Omega,
 \end{equation}
 where $a_k(z)$, $b(z)$ are smooth functions and
 \begin{equation}
 a_k(z)=-2\frac{<\nu,e_1>}{<\beta,\nu>}(\delta_i\beta_k)\tau_i(e_1),\ \ b(z)=2\frac{<\nu,e_1>}{<\beta,\nu>}\delta_i\varphi\tau_i(e_1),\
 \ \ on\ \  \partial\Omega.
 \end{equation}
So $V(0)=0$ on $\partial\Omega.$

Assume $Du(0)=0$, or we let
\begin{equation}H(z)=\frac{D_{11}u(z)-V(z)}{D_{11}u(0)}+\frac{2\nu_1\beta_1^T}{\beta\cdot\nu}+\widehat{B}|Du(z)-Du(0)|^2+\widehat{G}|z|^2-\widehat{A}x_n-1.\end{equation}

According to (\ref{2.3}) and (\ref{2.4}), we have
$LD_{11}u(z)\geq -\widetilde{C_5}$.
 \begin{equation}
 \begin{array}{ll}
 &-LV\\
 &=-L(a_kD_ku+b)\\
 &=-F^{i\bar j}\partial_i(\partial_{\bar j}a_kD_ku+a_k\partial_{\bar j}D_ku+\partial_{\bar j}b)\\
 &=-F^{i\bar j}[\partial_{i\bar j}a_kD_ku+\partial_{\bar j}a_k\partial_iD_ku\\
 & \ \ \ \ \ \ \ \        +\partial_ia_k\partial_{\bar j}D_ku+a_k\partial_{i\bar j}D_ku+
 \partial_{i\bar j}b]\\
 &\geq -F^{i\bar j}(\partial_{\bar j}a_k\partial_i D_ku+\partial_ia_k\partial_{\bar j}D_ku)
 -Ctr(F^{i\bar j}) -a_kD_k\widetilde{f}-C\\
 &\geq -\widetilde{C_6}tr(F^{i\bar j})-F^{i\bar j}\partial_iD_ku\partial_{\bar
 j}D_ku-\widetilde{C_7},
 \end{array}
 \end{equation}
 where $\widetilde{C_6}$ and $\widetilde{C_7}$ depend on $a_k$, $b$, $f^{\frac{1}{n}}$and
$|u|_{1,\overline\Omega}$.
  \begin{equation}
  \begin{array}{ll}
  L(|Du|^2)
  &=F^{i\bar j}(2\partial_iD_ku\partial_{\bar j}D_ku+2D_ku\partial_{i\bar j}D_ku)\\
  &=2F^{i\bar j}\partial _iD_ku\partial_{\bar j}D_ku+2D_k\widetilde{f}D_ku\\
  &\geq2F^{i\bar j}\partial _iD_ku\partial_{\bar j}D_ku-\widetilde{C_8},
  \end{array}
  \end{equation}
  where $\widetilde{C_8}$ depends on $|u|_{1,\overline\Omega}$ and $f^{\frac{1}{n}}$.
  \begin{equation}L\widehat{G}|z|^2=\widehat{G}tr(F^{i\bar j}).\end{equation}
  And see (2.26) in \cite{U2},
  \begin{equation}
  L(\frac{2\nu_1\beta^T}{\beta\cdot\nu})\geq-(\widetilde{C_{9}}\sqrt{\mu_2}|\beta|_{2,\overline{\Omega}}+
  \widetilde{C_{10}}|\beta^T|_{1,\overline{\Omega}})tr( F^{i\bar j})\ \ in \ S_{\mu_2}.\end{equation}
  Therefore,
\begin{equation}
\begin{array}{ll}
LH&\geq\frac{-\widetilde{C_6}tr(F^{i\bar j})-\widetilde{C_5}-\widetilde{C_7}-F^{i\bar j}\partial _iD_ku\partial_{\bar j}D_ku}{D_{11}u(0)}+2\widehat{B}F^{i\bar j}\partial _iD_ku\partial_{\bar j}D_ku-\widehat{B}\widetilde{C_8}+\widehat{G}tr(F^{i\bar j})\\
&-(\widetilde{C_{9}}\sqrt{\mu_2}|\beta|_{2,\overline{\Omega}}+
  \widetilde{C_{10}}|\beta^T|_{1,\overline{\Omega}})tr( F^{i\bar j})\\
&=(\frac{\widehat{G}}{2}tr(F^{i\bar j})-\widehat{B}\widetilde{C_8}-\frac{\widetilde{C_5}+\widetilde{C_7}}{D_{11}u(0)})\\&+(\frac{\widetilde{G}}{2}
-\frac{\widetilde{C_6}}{D_{11}u(0)}-\widetilde{C_{9}}\sqrt{\mu_2}|\beta|_{2,\overline{\Omega}}-
  \widetilde{C_{10}}|\beta^T|_{1,\overline{\Omega}})tr(F^{i\bar j})\\&+(2\widehat{B}-\frac{1}{D_{11}u(0)})F^{i\bar j}\partial _iD_ku\partial_{\bar j}D_ku\\
&\geq (\frac{\widehat{G}}{2}-\widehat{B}\widetilde{C_8}-\frac{\widetilde{C_5}+\widetilde{C_7}}{D_{11}u(0)})+(\frac{\widetilde{G}}{2}
-\frac{\widetilde{C_6}}{D_{11}u(0)}-\widetilde{C_{9}}\sqrt{\mu_2}|\beta|_{2,\overline{\Omega}}-
  \widetilde{C_{10}}|\beta^T|_{1,\overline{\Omega}})tr(F^{i\bar j})\\
&+(2\widehat{B}-\frac{1}{D_{11}u(0)})F^{i\bar j}\partial _iD_ku\partial_{\bar j}D_ku,\ in \ S_{\mu_2}.
 \end{array}
 \end{equation}
So if
\begin{equation}\label{cs1}
\frac{\widehat{G}}{2}-\widehat{B}\widetilde{C_8}-\frac{\widetilde{C_5}+\widetilde{C_7}}{D_{11}u(0)}\geq0,\ \frac{\widetilde{G}}{2}
-\frac{\widetilde{C_6}}{D_{11}u(0)}-\widetilde{C_{9}}\sqrt{\mu_2}|\beta|_{2,\overline{\Omega}}-\widetilde{C_{10}}|\beta^T|_{1,\overline{\Omega}}\geq0,\ 2\widehat{B}\geq\frac{1}{D_{11}u(0)},
\end{equation}
we have $LH\geq0$ in $S_{\mu_2}$.

On $\partial\Omega\cap\overline{S_{\mu_2}}$ near $0$,
from (\ref{4.19}), Lemma 4.3 and Lemma 4.4,
\begin{equation}
\begin{array}{ll}
H(z)&\leq 1+\frac{\widetilde{C_{12}}|z|^2}{D_{11}u(0)}+(\widehat{G}-\widehat{A}a)|z|^2-1\\
&\leq[\frac{\widetilde{C_{12}}}{D_{11}u(0)}+\widehat{G}-\widehat{A}a]|z|^2,
\end{array}
\end{equation}
if \begin{equation}\label{cs2}\frac{\widetilde{C_{12}}}{D_{11}u(0)}+\widehat{G}\leq\widehat{A}a,\end{equation} we have $H(z)\leq 0$ on $\partial\Omega\cap\overline{S_\mu}$ near $0$.

On the other hand, from (\ref{3.6}), (\ref{4.3}), and Theorem \ref{ge},
\begin{equation}
\begin{array}{ll}
H(z)&\leq (1+\frac{2}{\beta_0}|\beta^T|)+\frac{\widetilde{C_{13}}}{D_{11}u(0)}+\widetilde{C_{14}}\sqrt{\mu_2}|\beta^T|+\widehat{B}C+\widehat{G}|z|^2-\widehat{A}\mu_2-1,\\
&\leq \frac{2}{\beta_0}|\beta^T|+\frac{\widetilde{C_{13}}}{D_{11}u(0)}+\widetilde{C_{14}}\sqrt{\mu_2}|\beta^T|+\widehat{G}\widetilde{C_{16}}\mu_2-\widehat{A}\mu_2.
\end{array}
\end{equation}
if \begin{equation}\label{cs3}\frac{2}{\beta_0}|\beta^T|+\frac{\widetilde{C_{13}}}{D_{11}u(0)}+\widetilde{C_{14}}\sqrt{\mu_2}|\beta^T|+\widehat{B}C
+\widehat{G}\widetilde{C_{16}}\mu_2\leq \widehat{A}\mu_2,\end{equation} we have
$H(z)\leq 0$ on $\{x_n=\mu\}\cap \overline{S_{\mu_2}}$.

We now proceed to fix $\widehat{G}$, $\mu_2$ and $\widehat{B}$, depending on $\widehat{A}$ which will be fixed later. We first fix $\widehat{G}>0$ so small that
\begin{equation}
\widehat{G}\leq\frac{\widehat{A}a}{2}\ \ and \ \ \widehat{G}\widehat{C_{16}}\leq \frac{\widehat{A}}{2},
\end{equation}
and then fix $\mu_2\in(0,1)$ so that
\begin{equation}
 \widetilde{C_9}\sqrt{\mu_2}|\beta|_{2,\overline{\Omega}}\leq\frac{\widehat{G}}{4},
\end{equation}
and then fix $\widehat{B}$ so that
\begin{equation}
\widehat{B}\widetilde{C_8}\leq\frac{\widehat{G}}{4}\ \ and\ \widehat{B}C\leq\frac{\widehat{A}\mu_2}{4}.
\end{equation}
Then (\ref{cs1}), (\ref{cs2}) and (\ref{cs3}) will be hold whenever
\begin{equation}\label{cs4}
\begin{array}{ll}
&\frac{\widetilde{C_5}+\widetilde{C_7}}{D_{11}u(0)}\leq\frac{\widetilde{G}}{4},\\
&\frac{\widetilde{C_9}}{D_{11}u(0)}+\widetilde{C_{10}}|\beta^T|_{1,\overline{\Omega}}\leq\frac{\widehat{G}}{4},\\
&\frac{1}{D_{11}u(0)}\leq 2\widehat{B},\\
&\frac{\widetilde{C_{12}}}{D_{11}u(0)}\leq \frac{\widehat{A}a}{2},\\
&(\frac{2}{\beta_0}+\widetilde{C_{14}})|\beta^T|+\frac{\widetilde{C_{13}}}{D_{11}u(0)}\leq \frac{\widehat{A}\mu_2}{4}.
\end{array}
\end{equation}
When these conditions (\ref{cs4}) are satisfied we have $Q\leq 0$ in $S_{\mu_2}$ by the maximum principle, and since $H(0)=0$, then
$D_\beta H(0)\geq 0$.

 We have
\begin{equation}\label{q1}
  \begin{array}{ll}
D_\beta D_{11}u(0)&\leq D_{11}\varphi(0,u(0))-\sum_{k=1}^{2n}(D_{11}\beta_k)
  D_ku(0)-2\sum_{k=1}^{2n}(D_{1}\beta_k)D_kD_1u(0)+C\\
  &\leq\varphi_uD_{11}u(0)-\sum_{k=1}^{2n}(D_{11}\beta_k)
  D_ku(0)-2\sum_{k=1}^{2n}(D_{1}\beta_k)D_kD_1u(0)+C\\
  &\leq -\gamma_0D_{11}u(0)-2\sum_{k=1}^{2n}(D_1\beta_k)D_kD_1u(0)+C\\
  &=C-\gamma_0D_{11}u(0)-2\sum_{k=1}^{2n}(D_1\nu_k)D_kD_1u(0)-2\sum_{k=1}^{2n}D_1(\beta_k-\nu_k)D_kD_1u(0)\\
  &=C-\gamma_0D_{11}u(0)-2\sum_{k=1}^{2n-1}\frac{\partial^2 r}{\partial t_k\partial
  t_1}D_kD_1u(0)\\&-2\sum_{k=1}^{2n-1}D_1(\beta_k-\nu_k)D_kD_1u(0)
  -2[\frac{\partial^2 r}{\partial t_{2n}\partial
  t_1}+D_1(\beta_{2n}-\nu_{2n})]D_\nu D_1u(0).
  \end{array}
  \end{equation}
  To handle the last term above, we express $\nu$ in terms of
  $\beta$ and tangential components,
  \begin{equation}\label{q2}
  D_{\nu}D_1u(0)=\sum_{k=1}^{2n-1}(-\frac{\beta_k}{\beta_{2n}})D_kD_1u(0)+\frac{1}{\beta_{2n}}D_\beta D_1u(0).
  \end{equation}
  We take (\ref{q2}) into the last term of (\ref{q1}),
\begin{equation}\label{a}
  \begin{array}{ll}
  D_\beta D_{11}u(0)&\leq C-\gamma_0D_{11}u(0)-2\sum_{k=1}^{2n-1}\frac{\partial^2 r}{\partial t_k\partial
  t_1}D_kD_1u(0)-2\sum_{k=1}^{2n-1}D_1(\beta_k-\nu_k)D_kD_1u(0)\\
  &-2[\frac{\partial^2 r}{\partial t_{2n}\partial
  t_1}+D_1(\beta_{2n}-\nu_{2n})][\sum_{k=1}^{2n-1}(-\frac{\beta_k}{\beta_{2n}})D_{k1}u(0)+\frac{1}{\beta_{2n}}D_{1\beta}u(0)]\\
  &\leq C-\gamma_0D_{11}u(0)-2[\sum_{k=1}^{2n-1}\frac{\partial^2 r}{\partial t_k\partial
  t_1}+\sum_{k=1}^{2n-1}D_1(\beta_k-\nu_k)\\
  &+\sum_{k=1}^{2n-1}(\frac{\partial^2 r}{\partial t_{2n}\partial
  t_1}+D_1(\beta_{2n}-\nu_{2n}))(-\frac{\beta_k}{\beta_{2n}})]D_kD_1u(0)\\
  &=C-\gamma_0D_{11}u(0)-6\Lambda
  D_{11}u(0)+\sum_{k=1}^{2n-1}2[\Lambda\delta_{k1}-\frac{\partial^2 r}{\partial t_k\partial
  t_1}]D_kD_1u(0)\\
  &+\sum_{k=1}^{2n-1}2[\Lambda\delta_{k1}-D_1(\beta_k-\nu_k)]D_kD_1u(0)\\
  &+\sum_{k=1}^{2n-1}2[\Lambda\delta_{k1}+
  (\frac{\partial^2 r}{\partial t_{2n}\partial
  t_1}+D_1(\beta_{2n}-\nu_{2n}))(\frac{\beta_k}{\beta_{2n}})]D_kD_1u(0).
  \end{array}
  \end{equation}
Let \begin{equation}
A_1'=\sum_{k=1}^{2n-1}2[\Lambda\delta_{k1}-\frac{\partial^2
r}{\partial t_k\partial
  t_1}]D_kD_1u(0),
\end{equation}
\begin{equation}
A_2'=\sum_{k=1}^{2n-1}2[\Lambda\delta_{k1}-D_1(\beta_k-\nu_k)]D_kD_1u(0),
\end{equation}
\begin{equation}
A_3'=\sum_{k=1}^{2n-1}2[\Lambda\delta_{k1}+
  (\frac{\partial^2 r}{\partial t_{2n}\partial
  t_1}+D_1(\beta_{2n}-\nu_{2n}))(\frac{\beta_k}{\beta_{2n}})]D_kD_1u(0).
\end{equation}
By using an argument similar to that given in the
  proof of Theorem \ref{ge}, we want to obtain the
  relationship between $A_1',\ A_2' $ and $A_3'$ with $D_{11}u(0)$.

First, by the conclusion in \cite{L1} directly, we have
\begin{equation}\label{a1}
\begin{array}{ll}
A_1'&=(2e'_1[\Lambda I_{2n-1}-(\frac{\partial^2 r(0)}{\partial t_k\partial t_l})^{2n-1}_{k,l=1}])\cdot(D_1D_1 u(0),\cdots,D_{2n-1}D_1u(0))\\
&=(2e'_1 UU^t[\Lambda I_{2n-1}-(\frac{\partial^2 r(0)}{\partial t_k\partial t_l})^{2n-1}_{k,l=1}]U)\cdot((D_1D_1 u(0),\cdots,D_{2n-1}D_1u(0))U)\\
&\leq 2(\Lambda-\lambda_2(0))D_{11}u(0)+C,
\end{array}
\end{equation}
where $U$ is defined in (\ref{U}).

Because the orthogonal transformation does not change the distance between two points, the second identity of (\ref{a1}) holds.
One can find the detail for the inequality of (\ref{a1}) in \cite{L}( Theorem 3.5 and Theorem 4.2).

An argument similar to that given in the proof of Theorem \ref{ge}
gives,
\begin{equation}\label{a2'}
A_2'\leq 2(\Lambda+\frac{1+\sqrt{8n-7}}{2}\epsilon_0)D_{11}u(0)+C,
\end{equation}when $\epsilon_0\leq\frac{\Lambda}{1+\sqrt{8n-7}}$.
  And, if we take $\epsilon_0\leq\frac{\Lambda|\beta_{2n}|}{M(1+\sqrt{8n-7})}$, we have
  \begin{equation}\label{a3'}
  \begin{array}{ll}
 A_3'\leq
 2(\Lambda+|\frac{M}{\beta_{2n}}|\frac{1+\sqrt{8n-7}}{2}\epsilon_0)D_{11}u(0)+C,
  \end{array}\end{equation}
  where
  \begin{equation}M(0)=\frac{\partial^2r}{\partial t_{2n}\partial
  t_1}+D_1(\beta_{2n}-\nu_{2n}),\end{equation} and $|M(0)|\leq M$,
 $M$ is a constant independent of $0$.

Thus by inserting (\ref{a1}), (\ref{a2'}) and (\ref{a3'}) into
(\ref{a}) we get
\begin{equation}
\begin{array}{ll}
 D_\beta D_{11}u(0)&\leq C-\gamma_0D_{11}u(0) -6\Lambda
  D_{11}u(0)+2(\Lambda-\lambda_2(0))D_{11}u(0)\\
  &+2(\Lambda+\frac{1+\sqrt{8n-7}}{2}\epsilon_0)D_{11}u(0)
+2(\Lambda+|\frac{M}{\beta_{2n}}|\frac{1+\sqrt{8n-7}}{2}\epsilon_0)D_{11}u(0)\\
&=C-[\gamma_0+2\lambda_2(0)-(1+\sqrt{8n-7})\epsilon_0-(1+\sqrt{8n-7})|\frac{M}{\beta_{2n}}|\epsilon_0]D_{11}u(0).
  \end{array}
\end{equation}
Again, we take $\epsilon_0\leq
\frac{\gamma_0+2\lambda_2(0)}{2[(1+\sqrt{8n-7})(1+|\frac{M}{\beta_{2n}}|)]}$,
then
\begin{equation}
\frac{\gamma_0+2\lambda_2(z_0)}{2}\geq
(1+\sqrt{8n-7})\epsilon_0+(1+\sqrt{8n-7})|\frac{M}{\beta_{2n}}|\epsilon_0.
\end{equation}
Ultimately, we choose \begin{equation}\begin{array}{ll}
\epsilon_0&\leq
\min\{\frac{\Lambda}{1+\sqrt{8n-7}},\frac{\Lambda\beta_0}{M(1+\sqrt{8n-7})},\frac{\gamma_0+2\lambda_2(z_0)}{2[(1+\sqrt{8n-7})(1+\frac{M}{\beta_0})]}\}\\
&\leq\min\{\frac{\Lambda}{1+\sqrt{8n-7}},\frac{\Lambda|\beta_{2n}|}{M(1+\sqrt{8n-7})},\frac{\gamma_0+2\lambda_2(z_0)}{2[(1+\sqrt{8n-7})(1+|\frac{M}{\beta_{2n}}|)]}\},
\end{array}
\end{equation}
then
\begin{equation} \gamma_0+2\lambda_2(0)-(1+\sqrt{8n-7})\epsilon_0-(1+\sqrt{8n-7})|\frac{M}{\beta_{2n}}|\epsilon_0\geq\frac{\sigma_1}{2},
\end{equation}
where $\sigma_1=\inf_{\partial\Omega}(\gamma_0+2\lambda_2)$.
Then
\begin{equation}
\begin{array}{ll}
0&\leq D_\beta H(0)\\
&\leq \frac{C-\frac{\sigma_1}{2}D_{11}u(0)}{D_{11}u(0)}-\widehat{A}D_\beta x_n(0).
\end{array}
\end{equation}
\begin{equation}
[\frac{\sigma_1}{2}-\widehat{A}\beta\cdot\nu(0)]D_{11}u(0)\leq C.
\end{equation}
Finally we fix $\widehat{A}$ so small that
\begin{equation}
\widehat{A}\beta\cdot\nu\leq\frac{\sigma_1}{4}, on\ \partial\Omega,
\end{equation}
 we have for any tangential vector  $D_{11}u(0)\leq \frac{\sigma_1}{4}.$ We finish the proof of Theorem 4.1.

\begin{remark}
We remark that the condition (\ref{ob2'}) is necessary in the above
proof. It can not be extended to the degenerate oblique case
$(\beta\cdot\nu)\geq0$.
\end{remark}
\begin{remark}
If $\gamma_0$ is sufficiently large, then we do not need any
structural assumptions on $\beta$ and the principal curvature of
$\partial \Omega$, for example, the condition (\ref{ob2}) can be
omitted.
\end{remark}

\begin{theorem}\label{thm4}Let $\Omega$ be a bounded strictly pseudoconvex domain in
$\bf C^n$ with $C^4$ boundary. Let $\varphi\in C^{3,1}(\partial
\Omega\times R)$, and $f\in C^2(\overline{\Omega})$, so that they
satisfy (\ref{ob2'})-(\ref{ob6}). In addition $\varphi$ satisfies
(\ref{ob3}). Then
\begin{equation}\label{3.12}|u|_{2,\overline{\Omega}}\leq C,\end{equation} where $C$ is a
constant depending on $\
|f^{\frac{1}{n}}|_{C^{2}(\overline\Omega)},\ \gamma_0,\ \Omega,
\beta \ and \ |\varphi|_{C^{3,1}(\partial\Omega)}$.\end{theorem}
\section{The Proof of Theorem}

In this subsection, we shall prove Theorem\ref{thm}. Although, the
argument in this part is rather standard, we present its sketch here
for completeness.

With the $C^0$, $C^1$ and $C^2$ bounds for the solution $u$
formulated in previous sections, the complex Monge-Amp\`ere equation
(\ref{cma}) is uniformly elliptic. Since the second derivatives are
bounded, the bounds of their H{\"o}lder norms follow from the
uniformly elliptic theory developed by Lieberman and Trudinger
\cite{LieTru1986}, that is $\|u\|_{C^{2,\alpha}(\bar \Omega)}\leq C$
for some $\alpha\in (0,1)$. Such a {\it priori} $C^{2,\alpha}$
estimation enables us to carry out the method of continuity in
\cite{GT} and \cite{LTU}, thus we obtain the existence of classical
solutions. This completes the proof of Theorem \ref{thm}.
\begin{remark}
We remark that the proof for the uniqueness of the solutions is
similar to that in \cite{L1} for the Neumann boundary case, here we
need to use the condition $\varphi_u<0$.
\end{remark}

\bibliographystyle{amsplain}

\end{document}